# A COMPANION FOR THE KIEFER–WOLFOWITZ–BLUM STOCHASTIC APPROXIMATION ALGORITHM

By Abdelkader Mokkadem and Mariane Pelletier

*Université de Versailles–Saint-Quentin*

A stochastic algorithm for the recursive approximation of the location $\theta$ of a maximum of a regression function was introduced by Kiefer and Wolfowitz [*Ann. Math. Statist.* **23** (1952) 462–466] in the univariate framework, and by Blum [*Ann. Math. Statist.* **25** (1954) 737–744] in the multivariate case. The aim of this paper is to provide a companion algorithm to the Kiefer–Wolfowitz–Blum algorithm, which allows one to simultaneously recursively approximate the size $\mu$ of the maximum of the regression function. A precise study of the joint weak convergence rate of both algorithms is given; it turns out that, unlike the location of the maximum, the size of the maximum can be approximated by an algorithm which converges at the parametric rate. Moreover, averaging leads to an asymptotically efficient algorithm for the approximation of the couple $(\theta, \mu)$.

**1. Introduction.** Consider two random variables $X$ and $Z$ with values in $\mathbb{R}^d$ and $\mathbb{R}$, respectively, that have unknown common distribution $\mathbb{P}_{X,Z}$. Assume that the regression function $f(\cdot) = \mathbb{E}(Z|X = \cdot): \mathbb{R}^d \to \mathbb{R}$ exists, is sufficiently smooth and has a unique maximizer $\theta \in \mathbb{R}^d$,

$$\theta = \arg\max_{x \in \mathbb{R}^d} \mathbb{E}(Z|X = x),$$

and assume that observations $Z(x)$ of $f(x)$ are available at any level $x$ [$Z(x)$ has conditional distribution $\mathcal{L}(Z|X = x)$]. Kiefer and Wolfowitz [15] (in the case $d = 1$) and Blum [1] (in the case $d \geq 1$) have introduced an algorithm, which allows one to recursively approximate $\theta$. Their procedure consists in running the recursion

$$(1) \qquad \theta_{n+1} = \theta_n + a_n Y_n,$$

---

where $(a_n)$ is a positive nonrandom sequence that goes to zero as $n$ goes to infinity, and $Y_n$ is a (random) approximation of $\nabla f(\theta_n)$, the gradient of $f$ at the point $\theta_n$. More precisely, let $(c_n)$ be a positive nonrandom sequence that goes to zero, and let $(e_1, \ldots, e_d)$ denote the canonical basis of $\mathbb{R}^d$; the approximation $Y_n$ introduced by Kiefer and Wolfowitz [15] and Blum [1] is the $d$-dimensional vector

$$Y_n = \frac{1}{2c_n} \{ Z(\theta_n + c_n e_i) - Z(\theta_n - c_n e_i) \}_{i \in \{1, \ldots, d\}}.$$

Kiefer and Wolfowitz [15] proved the convergence in probability of $\theta_n$ to $\theta$ and Blum [1] established its almost sure convergence. Their algorithm (1) has since been widely studied and their pioneering work extended in many directions. Among many, let us cite Fabian [11], Kushner and Clark [17], Hall and Heyde [13], Ruppert [31], Chen [3], Spall [33, 34], Polyak and Tsybakov [30], Dippon and Renz [8], Pelletier [26], Chen, Duncan and Pasik-Duncan [4] and Dippon [6].

As noted by Kiefer and Wolfowitz [15], the statistical importance of approximating the maximizer $\theta$ of the regression function $f$ is obvious and need not be discussed. Although the approximation of the size of the maximum, that is of the parameter $\mu = f(\theta)$, seems important as well, this problem has, as far as we know, never been considered. The aim of this paper is to propose an algorithm, which by using the approximation $\theta_n$ of $\theta$ defined by (1), allows one to simultaneously recursively approximate $\mu$ by a sequence $\mu_n$ that converges almost surely to $\mu$, and to study the joint weak convergence rate of $\theta_n$ and $\mu_n$.

The algorithm we present to approximate $\mu$ is defined by

$$(2) \qquad \mu_{n+1} = (1 - \tilde{a}_n) \mu_n + \tilde{a}_n \tilde{Y}_n,$$

where $(\tilde{a}_n)$ is a positive nonrandom sequence that goes to zero as $n$ goes to infinity, and $\tilde{Y}_n$ is an approximation of $f(\theta_n)$. This approximation method has certain similarities to the sequential procedure for estimating discontinuities of a regression function or surface proposed by Hall and Molchanov [14]. A first way to approximate $f(\theta_n)$ is to take the average of the observations of $f(\theta_n + c_n e_i)$ and $f(\theta_n - c_n e_i)$ used for the computation of $Y_n$; all these observations or only a symmetric part of them may be used. More precisely, let $\mathcal{S}$ denote a (nonempty) subset of $\{1, 2, \ldots, d\}$, and define the real-valued random sequence $(\tilde{Y}_n)$ by

$$\tilde{Y}_n = \frac{1}{\delta} \sum_{i \in \mathcal{S}} \{ Z(\theta_n + c_n e_i) + Z(\theta_n - c_n e_i) \},$$

where $\delta$ is twice the number of elements in $\mathcal{S}$. Note that in the case the step size in (2) is chosen such that $(\tilde{a}_n) \equiv (n^{-1})$ and if $\mathcal{S} = \{1, 2, \ldots, d\}$, then $\mu_{n+1}$ is simply the average of all the observations made for the approximation $\theta_n$ of $\theta$, that is,

$$\mu_{n+1} = \frac{1}{n}\sum_{k=1}^{n}\tilde{Y}_k = \frac{1}{2dn}\sum_{i \in \{1,2,\dots,d\}, k \in \{1,2,\dots,n\}}\{Z(\theta_k + c_k e_i) + Z(\theta_k - c_k e_i)\}.$$

We prove that, under suitable assumptions, $\mu_n$ converges almost surely to $\mu$. Moreover, we study the weak convergence rate of the couple $(\theta_n - \theta, \mu_n - \mu)$. As was already well known, the optimal convergence rate of $\theta_n$ (which is $n^{1/3}$) is obtained by choosing in (1) $(a_n) \equiv (a_0 n^{-1})$ with adequate conditions on $a_0$, and $(c_n) \equiv (c_0 n^{-1/6})$, $c_0 > 0$; setting $(\tilde{a}_n) \equiv (\tilde{a}_0 n^{-1})$, $\tilde{a}_0 > 1/2$, in (2) then makes $\mu_n$ converge with the rate $n^{1/3}$ also. Now, other choices of $(c_n)$ in (1) and (2) allow one to obtain a convergence rate for $\mu_n$ close to (but less than) the parametric rate $\sqrt{n}$; however, in this case, the convergence rate of $\theta_n$ becomes close to $n^{1/4}$. This constatation makes clear the drawback of the double algorithm (1) and (2): when choosing the sequence $(c_n)$ [or, in other words, the points where the observations $Z(\theta_k \pm c_k e_i)$ of $f(\theta_k \pm c_k e_i)$ are taken], a compromise must be made since both sequences $\theta_n$ and $\mu_n$ cannot simultaneously converge at the optimal rate.

The way to address this drawback is of course not to use the same sequence $(c_n)$ (i.e., to use different observations) for the approximation of $\nabla f(\theta_n)$ in (1) on the one hand, and for the approximation of $f(\theta_n)$ in (2) on the other hand. More precisely, let $\delta \geq 1$, $Z_i(\theta_n)$, $1 \leq i \leq \delta$, be $\delta$ independent observations of $f(\theta_n)$, $\tilde{\mathcal{Y}}_n$ be the approximation of $f(\theta_n)$ defined by

$$(3) \qquad \tilde{\mathcal{Y}}_n = \frac{1}{\delta}\sum_{i=1}^{\delta} Z_i(\theta_n),$$

and let the approximation algorithm for $\mu$ be defined as

$$(4) \qquad \mu_{n+1} = (1 - \tilde{a}_n)\mu_n + \tilde{a}_n \tilde{\mathcal{Y}}_n.$$

We prove that the sequence $\mu_n$ defined in this way still converges almost surely to $\mu$. Moreover, we study the joint weak convergence rate of $\theta_n$ and $\mu_n$ defined by (1) and (4), respectively. We prove in particular that if the stepsizes in (1) and (4) are chosen such that $(a_n) \equiv (a_0 n^{-1})$, with adequate conditions on $a_0$, $(c_n) \equiv (c_0 n^{-1/6})$, $c_0 > 0$, and $(\tilde{a}_n) \equiv (\tilde{a}_0 n^{-1})$, $\tilde{a}_0 > 1/2$, then $(\theta_n)$ converges with its optimal rate $n^{1/3}$, and $(\mu_n)$ with the parametric rate $\sqrt{n}$. Moreover, choosing $\tilde{a}_0 = 1$ leads to the minimum asymptotic variance of $(\mu_n)$: when $(\tilde{a}_n) \equiv (n^{-1})$, the algorithm (4) is asymptotically efficient. Note that this case corresponds to the case

$$\mu_{n+1} = \frac{1}{n}\sum_{k=1}^{n}\tilde{\mathcal{Y}}_k = \frac{1}{\delta n}\sum_{i \in \{1,2,\dots,\delta\}, k \in \{1,2,\dots,n\}} Z_i(\theta_k).$$

The striking aspect of our result on (4) is that, whereas approximation of the size of the maximum of a regression function is typically a nonparametric problem, and although the stochastic approximation algorithm (4) uses approximation of the location of the maximum of the regression function $\theta_n$ (which itself does not converge with the parametric rate), the convergence rate we obtain for the sequence $\mu_n$ is the parametric rate $\sqrt{n}$. This is explained by the fact that although $\mu_n$ depends (through $\hat{\mathcal{Y}}_n$) on $\theta_n$, the quantity which actually is involved in the convergence rate of $(\mu_n)$ is $\|\theta_n - \theta\|^2$, and, for suitable choices of $(a_n)$ and $(c_n)$, this quantity goes to zero faster than $\sqrt{n}$. [Of course, this is still true in the framework of the double algorithm (1) and (2), but in this case the convergence rate of $(\mu_n)$ depends on $(c_n)$ and is less than $\sqrt{n}$.]

Now, as is well known, the choice of the step size $(a_n) \equiv (a_0 n^{-1})$ in (1) is the one which leads to the optimal convergence rate of $\theta_n$, but it induces conditions on $a_0$ which are difficult to handle because of dependence on an unknown parameter [see (9) in the sequel]. The well known approach used to obtain optimal convergence rates for stochastic approximation algorithms without a tedious condition on the step size is to use the averaging principle independently introduced by Ruppert [32] and Polyak [28]. Their averaging procedure, which has been widely discussed and extended (see, among many others, Yin [35], Delyon and Juditsky [5], Polyak and Juditsky [29], Kushner and Yang [18], Le Breton [19], Le Breton and Novikov [20], Dippon and Renz [7, 8] and Pelletier [27]), allows one to obtain asymptotically efficient algorithms, that is, algorithms which not only converge at the optimal rate, but also have an optimal asymptotic covariance matrix. This procedure consists in (i) running the approximation algorithm by using slower step sizes and (ii) computing a suitable average of the approximations obtained in (i).

Let us now give our scheme to efficiently approximate $\theta$ and $\mu$ simultaneously. First, we apply the averaging principle to the approximating algorithm (1) of $\theta$ by proceeding as follows. Let the step size $(a_n)$ in (1) satisfy $\lim_{n \to \infty} na_n = \infty$, let the sequence $(\theta_k)$ be defined by the algorithm (1) and set

$$(5) \qquad \overline{\theta}_n = \frac{1}{\sum_{k=1}^n c_k^2} \sum_{k=1}^n c_k^2 \theta_k.$$

It is well known that the sequence $(\overline{\theta}_n)$ is asymptotically efficient (see, e.g., [8]). Then, to approximate $\mu$ efficiently, we can just set $(\tilde{a}_n) \equiv (n^{-1})$ in (4) since this algorithm is asymptotically efficient (see the comments below Theorem 2). However, when adding observations of $f$, it seems more natural to take the observations at the point $\overline{\theta}_n$ (rather than at $\theta_n$) since $\overline{\theta}_n$ converges to $\theta$ faster than $\theta_n$ does. That is the reason why we let $\delta \geq 1$, $Z_i(\overline{\theta}_n)$, $1 \leq$ $i \leq \delta$, be $\delta$ independent observations of $f(\overline{\theta}_n)$, $\overline{\mathcal{Y}}_n$ be the approximation of $f(\overline{\theta}_n)$ defined by

$$(6) \qquad \overline{\mathcal{Y}}_n = \frac{1}{\delta} \sum_{i=1}^{\delta} Z_i(\overline{\theta}_n),$$

and let the approximation algorithm for $\mu$ be defined as

$$(7) \qquad \mu_{n+1} = \left(1 - \frac{1}{n}\right)\mu_n + \frac{1}{n}\overline{\mathcal{Y}}_n.$$

The consistency of $\mu_n$ defined by (7) is obvious; we study the joint weak asymptotic behavior of $\overline{\theta}_n$ and $\mu_n$ defined by (5) and (7). We prove in particular that by setting $(c_n) \equiv (c_0 n^{-1/6})$ in (1), we obtain simultaneously the asymptotic efficiency of both sequences $(\overline{\theta}_n)$ and $(\mu_n)$.

Let us finally mention that, in the case where no additional observations are taken to approximate $\mu$, we can of course also average the algorithm (1). However, we shall point out that when the only parameter of interest in the double algorithm (1) and (2) is $\mu$, it is preferable not to do so. As a matter of fact, we show there are possible choices of $(a_n)$ for which there is no tedious condition on $a_0$, and which lead to better convergence rates for $(\mu_n)$ than those which can be reached by averaging $\theta_n$.

## 2. Assumptions and main results.
Let us first define the class of positive sequences that will be used in the statement of our assumptions.

DEFINITION 1. Let $\alpha \in \mathbb{R}$ and $(v_n)$ be a nonrandom positive sequence. We say that $(v_n) \in \mathcal{GS}(\alpha)$ if

$$(8) \qquad \lim_{n \to \infty} n\left[1 - \frac{v_{n-1}}{v_n}\right] = \alpha.$$

Condition (8) was introduced by Galambos and Seneta [12] to define regularly varying sequences (see also [2]). Typical sequences in $\mathcal{GS}(\alpha)$ are, for $a \in \mathbb{R}$, $n^\alpha(\log n)^a$, $n^\alpha(\log \log n)^a$, and so on.

Set

$$W_{n,i}^+ = Z(\theta_n + c_n e_i) - f(\theta_n + c_n e_i),$$

$$W_{n,i}^- = Z(\theta_n - c_n e_i) - f(\theta_n - c_n e_i),$$

$$\mathcal{W}_{n,i} = Z_i(\theta_n) - f(\theta_n),$$

$$\overline{\mathcal{W}}_{n,i} = Z_i(\overline{\theta}_n) - f(\overline{\theta}_n).$$

(The notation $\mathcal{W}_{n,i}$ (resp. $\overline{\mathcal{W}}_{n,i}$) is useful only in the case $(\mu_n)$ is defined by (4) [resp. by (7)].) In order to state our assumptions in a compact way, we introduce the sequence $(b_n)$ defined as

$$(b_n) \equiv \begin{cases} (c_n), & \text{in the case } (\mu_n) \text{ is defined by (2),} \\ 0, & \text{in the case } (\mu_n) \text{ is defined by (4) or by (7),} \end{cases}$$

and set

$$\mathcal{U}_{n,i} = \begin{cases} 0, & \text{in the case } (\mu_n) \text{ is defined by (2),} \\ \mathcal{W}_{n,i}, & \text{in the case } (\mu_n) \text{ is defined by (4),} \\ \overline{\mathcal{W}_{n,i}}, & \text{in the case } (\mu_n) \text{ is defined by (7).} \end{cases}$$

The assumptions to which we shall refer in the sequel are the following.

(A1) $\lim_{n \to \infty} \theta_n = \theta$ a.s.

(A2) $f$ is three-times continuously differentiable in a neighborhood of $\theta$, where the Hessian $D^2 f(\theta)$ of $f$ at $\theta$ is negative definite with maximal eigenvalue $-L^{(\theta)} < 0$.

(A3) Let $\mathcal{G}_n$ be the $\sigma$-field spanned by $\{W^+_{m,i}, W^-_{p,j}, \mathcal{U}_{q,k} \ 1 \le i,j \le d, \ 1 \le k \le \delta, \ 1 \le m,p,q \le n-1\}$.

    (i) $W^+_{n,i}$, $W^-_{n,j}$ and $\mathcal{U}_{n,k}$ $(i,j \in \{1,\dots,d\}, \ k \in \{1,\dots,\delta\})$ are independent conditionally on $\mathcal{G}_n$.

    (ii) For some $\sigma > 0$, $\mathrm{Var}(Z|X = x) = \sigma^2$ for all $x \in \mathbb{R}^d$, while, for some $m > 2$, $\sup_{x \in \mathbb{R}^d} \mathbb{E}(|Z|^m|X = x) < \infty$.

(A4) (i) There exists $\alpha \in ]\max\{1/2, 2/m\}, 1]$ such that $(a_n) \in \mathcal{GS}(-\alpha)$.

    (ii) There exists $\tau \in ]0, \alpha/2[$ such that $(c_n) \in \mathcal{GS}(-\tau)$.

    (iii) $\lim_{n \to \infty} n a_n \in ]\max\{\frac{1-2\tau}{2L^{(\theta)}}, \frac{2\tau}{L^{(\theta)}}\}, \infty]$.

    (iv) There exists $\tilde{\alpha} \in ]\max\{1/2, 2/m\}, 1]$ such that $(\tilde{a}_n) \in \mathcal{GS}(-\tilde{\alpha})$.

    (v) • In the case $\lim_{n \to \infty} \tilde{a}_n^{-1} b_n^4 = 0$, we have $\lim_{n \to \infty} \tilde{a}_n^{-1/2} a_n \times \log(\sum_{k=1}^n a_k)/c_n^2 = 0$ and $\lim_{n \to \infty} \tilde{a}_n^{-1} c_n^8 = 0$.

      • In the case $\lim_{n \to \infty} \tilde{a}_n^{-1} b_n^4 \in ]0, \infty]$, we have $\sum \tilde{a}_n b_n^4 < \infty$ and $\lim_{n \to \infty} a_n \log(\sum_{k=1}^n a_k)/c_n^4 = 0$.

    (vi) $\lim_{n \to \infty} n \tilde{a}_n \in ]\frac{1}{2}, \infty]$.

*Comments on the assumptions.* (1) Theorem 3 in [1] ensures that (A1) holds under (A2)–(A4) and the following additional conditions: (i) $\alpha + \tau > 1$ and $2(\alpha - \tau) > 1$; (ii) $D^2 f$ is bounded; (iii) $\forall \delta > 0$, $\sup_{\|x-\theta\| \ge \delta} f(x) < f(\theta)$; (iv) $\forall \varepsilon > 0$, $\exists \rho(\epsilon) > 0$ such that $\|x - \theta\| \ge \varepsilon \Rightarrow \|\nabla f(x)\| \ge \rho(\epsilon)$. Let us underline that the conditions (i) on $\alpha$ and $\tau$ are satisfied as soon as $\alpha \in ]5/6, 1]$ and $\tau \in [1/6, 1/4]$, which include the most interesting choices of step sizes, as we shall see later on. Let us also mention that similar conditions, but which are less restrictive on $\alpha$ and $\tau$, can be found in [22] and [13]. Another kind of conditions with particular emphasis on control theory applications is given in [9, 17, 21]. The approach in these three references is to associate the approximation algorithm (1) with a deterministic differential equation in terms of which conditions are given to ensure (A1).

(2) Assumptions (A4)(i)–(iii) are the conditions on the step sizes required to establish the weak convergence rate of $\theta_n$; assumptions (A4)(iv)–(vi) are the additional ones needed for the consistency and for the weak convergence rate of $\mu_n$.

(3) Condition (A4)(iii) [resp. (A4)(vi)] requires $a_n = O(n^{-1})$ [resp. $\tilde{a}_n = O(n^{-1})$] and, in the case $(a_n) \equiv (a_0 n^{-1})$ [resp. $(\tilde{a}_n) \equiv (\tilde{a}_0 n^{-1})$],

$$(9) \qquad a_0 > \max\left\{\frac{1-2\tau}{2L^{(\theta)}}; \frac{2\tau}{L^{(\theta)}}\right\},$$

(resp. $\tilde{a}_0 > 1/2$). Set $\log_1(n) = \log n$ and, for $j \geq 1$, $\log_{j+1}(n) = \log[\log_j(n)]$. Our conditions allow the use of the step size $(a_n) \equiv (a_0[\log_p(n)]^\alpha n^{-1})$ introduced by Koval and Schwabe [16]; this step size has the advantage to lead to convergence rates very close to the ones obtained by using $(a_0 n^{-1})$, without requiring the tedious condition (9) on $a_0$.

(4) Assumption (A4)(v) is in particular satisfied as soon as the following conditions hold:

- If $\lim_{n\to\infty} \tilde{a}_n^{-1} b_n^4 = 0$, then $\frac{\tilde{\alpha}}{8} < \tau < \frac{\alpha}{2} - \frac{\tilde{\alpha}}{4}$.
- If $\lim_{n\to\infty} \tilde{a}_n^{-1} b_n^4 \in ]0, \infty]$, then $\frac{1-\tilde{\alpha}}{4} < \tau < \frac{\alpha}{4}$.

Our first result is the following proposition, which states the consistency of $\mu_n$ in the case $\mu_n$ is defined either by (2) or by (4).

PROPOSITION 1. *Let $\mu_n$ be defined either by (2) or by (4), and assume (A1)–(A3) and (A4)(i)–(v) are satisfied. Then we have $\lim_{n\to\infty} \mu_n = \mu$ a.s.*

In order to state the weak convergence rate of $(\theta_n^T, \mu_n)^T$, we set

$$\xi^{(\theta)} = (1-2\tau)\lim_{n\to\infty}(na_n)^{-1},$$

$$\zeta^{(\theta)} = 4\tau\lim_{n\to\infty}(na_n)^{-1},$$

$$\xi^{(\mu)} = \lim_{n\to\infty}(n\tilde{a}_n)^{-1},$$

$$\zeta^{(\mu)} = 4\tau\lim_{n\to\infty}(n\tilde{a}_n)^{-1},$$

$$(10) \qquad \Sigma^{(\theta)} = -\frac{\sigma^2}{4}\left[D^2f(\theta) + \frac{\xi^{(\theta)}}{2}I_d\right]^{-1},$$

$$(11) \qquad \Delta^{(\theta)} = -\frac{1}{6}\left[D^2f(\theta) + \frac{\zeta^{(\theta)}}{2}I_d\right]^{-1}\left\{\frac{\partial^3 f}{\partial x_i^3}(\theta)\right\}_{1\leq i\leq d},$$

$$(12) \qquad \Sigma^{(\mu)} = \frac{\sigma^2}{\delta(2-\xi^{(\mu)})},$$

(13) $$\Delta^{(\mu)} = \frac{2}{(2 - \zeta^{(\mu)})\delta} \sum_{i \in \mathcal{S}} \frac{\partial^2 f}{\partial x_i^2}(\theta),$$

where $\sigma^2$ is defined in (A3), and where $I_d$ denotes the $d \times d$ identity matrix. Let us underline that assumption (A4) implies that $\xi^{(\theta)}, \zeta^{(\theta)} \in [0, 2L^{(\theta)}[$ and $\xi^{(\mu)}, \zeta^{(\mu)} \in [0, 2[$; the parameters $\Sigma^{(\theta)}$, $\Delta^{(\theta)}$, $\Sigma^{(\mu)}$ and $\Delta^{(\mu)}$ are thus well defined.

We now state the joint weak convergence rate of $\theta_n$ and $\mu_n$ in the case $\mu_n$ where is defined by the algorithm (2).

THEOREM 1. *Let $(\mu_n)$ be defined by (2), and assume that* (A1)–(A4) *hold.*

(1) *If $\lim_{n \to \infty} a_n^{-1} c_n^6 = \infty$ and if $\lim_{n \to \infty} \tilde{a}_n^{-1} c_n^4 = \infty$, then*

$$\begin{pmatrix} c_n^{-2}(\theta_n - \theta) \\ c_n^{-2}(\mu_n - \mu) \end{pmatrix} \xrightarrow{\mathbb{P}} \begin{pmatrix} \Delta^{(\theta)} \\ \Delta^{(\mu)} \end{pmatrix}.$$

(2) *If there exists $\gamma_1 \geq 0$ such that $\lim_{n \to \infty} a_n^{-1} c_n^6 = \gamma_1$ and if $\lim_{n \to \infty} \tilde{a}_n^{-1} c_n^4 = \infty$, then*

$$\begin{pmatrix} \sqrt{a_n^{-1} c_n^2}(\theta_n - \theta) \\ c_n^{-2}(\mu_n - \mu) \end{pmatrix} \xrightarrow{\mathcal{D}} \begin{pmatrix} \mathcal{Z} \\ \Delta^{(\mu)} \end{pmatrix},$$

*where $\mathcal{Z}$ is $\mathcal{N}(\sqrt{\gamma_1}\Delta^{(\theta)}, \Sigma^{(\theta)})$-distributed.*

(3) *If $\lim_{n \to \infty} a_n^{-1} c_n^6 = \infty$ and if there exists $\gamma_2 \geq 0$ such that $\lim_{n \to \infty} \tilde{a}_n^{-1} c_n^4 = \gamma_2$, then*

$$\begin{pmatrix} c_n^{-2}(\theta_n - \theta) \\ \sqrt{\tilde{a}_n^{-1}}(\mu_n - \mu) \end{pmatrix} \xrightarrow{\mathcal{D}} \begin{pmatrix} \Delta^{(\theta)} \\ \mathcal{Z}' \end{pmatrix},$$

*where $\mathcal{Z}'$ is $\mathcal{N}(\sqrt{\gamma_2}\Delta^{(\mu)}, \Sigma^{(\mu)})$-distributed.*

(4) *If there exist $\gamma_1 \geq 0$ and $\gamma_2 \geq 0$ such that $\lim_{n \to \infty} a_n^{-1} c_n^6 = \gamma_1$ and $\lim_{n \to \infty} \tilde{a}_n^{-1} c_n^4 = \gamma_2$, then*

$$\begin{pmatrix} \sqrt{a_n^{-1} c_n^2}(\theta_n - \theta) \\ \sqrt{\tilde{a}_n^{-1}}(\mu_n - \mu) \end{pmatrix} \xrightarrow{\mathcal{D}} \mathcal{N}\left( \begin{pmatrix} \sqrt{\gamma_1}\Delta^{(\theta)} \\ \sqrt{\gamma_2}\Delta^{(\mu)} \end{pmatrix}, \begin{pmatrix} \Sigma^{(\theta)} & 0 \\ 0 & \Sigma^{(\mu)} \end{pmatrix} \right).$$

*Comments on Theorem* 1.

(1) As is already well known, the optimal convergence rate of $(\theta_n)$ is obtained by choosing $(a_n) \equiv (a_0 n^{-1})$, $a_0$ satisfying (9) and $(c_n) \equiv (n^{-1/6})$. In this framework, the best convergence rate of $(\mu_n)$ is $n^{1/3}$; it is obtained in the following ways:

- either $(\tilde{a}_n)$ is chosen such that $\lim_{n \to \infty} \tilde{a}_n^{-1} n^{-2/3} = \infty$, the convergence rate of $(\mu_n)$ being then given by part (2) of Theorem 1,

- or $(\tilde{a}_n) \equiv (n^{-2/3})$, the convergence rate of $(\mu_n)$ being then given by part (4) of Theorem 1.

(2) The optimal convergence rate of $(\mu_n)$ is close to (but less than) $\sqrt{n/\log\log n}$. More precisely, let $(v_n) \in \mathcal{GS}(0)$ be such that $\lim_{n\to\infty} v_n = \infty$. For $(\mu_n)$ to converge with the rate $\sqrt{n/(v_n\log\log n)}$, one must choose $(a_n) \equiv (a_0 n^{-1})$, $a_0$ satisfying (9), and

- either $(\tilde{a}_n) \equiv (n^{-1})$ and $(c_n) \equiv (v_n^{1/4}[\log\log n]^{1/4} n^{-1/4})$, the convergence rate of $(\mu_n)$ being then given by part (2) of Theorem 1,

- or $(\tilde{a}_n) \equiv (n^{-1} v_n \log\log n)$ and $(c_n) = O(\tilde{a}_n^{1/4})$, the convergence rate of $(\mu_n)$ being then given by part (4) of Theorem 1.

In this framework, the best convergence rate of $(\theta_n)$ is $n^{1/4} v_n^{1/4}[\log\log n]^{1/4}$.

(3) The tedious condition (9) on $a_0$ can be avoided by choosing $(a_n) \equiv (n^{-1}\log_p n)$. The convergence rate of $(\mu_n)$ is then close to (but less than) $\sqrt{n/(\log_p n\log\log n)}$. More precisely, let $(v_n) \in \mathcal{GS}(0)$ be such that $\lim_{n\to\infty} v_n = \infty$. For $(\mu_n)$ to converge with the rate $\sqrt{n/(v_n\log_p n\log\log n)}$, one can choose

- either $(\tilde{a}_n) \equiv (n^{-1})$ and $(c_n) \equiv (v_n^{1/4}[\log_p n]^{1/4}[\log\log n]^{1/4} n^{-1/4})$, the convergence rate of $(\mu_n)$ being then given by part (2) of Theorem 1,

- or $(\tilde{a}_n) \equiv (n^{-1} v_n \log_p n\log\log n)$ and $(c_n) = O(\tilde{a}_n^{1/4})$, the convergence rate of $(\mu_n)$ being then given by part (4) of Theorem 1.

In this case, the best convergence rate of $(\theta_n)$ is $n^{1/4} v_n^{1/4}[\log\log n]^{1/4}[\log_p n]^{-1/4}$.

The double algorithm (1) and (2) has thus two disadvantages: (i) it is not possible to choose a sequence $(c_n)$ such that the convergence rates of $(\theta_n)$ and $(\mu_n)$ are simultaneously optimal; (ii) the sequence $(\mu_n)$ cannot converge at the parametric rate.

We now state the joint weak convergence rate of $\theta_n$ and $\mu_n$ in the case additional observations are made for the computation of $\mu_n$, that is, in the case $(\mu_n)$ is defined by (4).

THEOREM 2. *Let $(\mu_n)$ be defined by (4), and assume that* (A1)–(A4) *hold.*

(1) *If $\lim_{n\to\infty} a_n^{-1} c_n^6 = \infty$, then*

$$\begin{pmatrix} c_n^{-2}(\theta_n - \theta) \\ \sqrt{\tilde{a}_n^{-1}}(\mu_n - \mu) \end{pmatrix} \xrightarrow{\mathcal{D}} \begin{pmatrix} \Delta^{(\theta)} \\ \mathcal{Z}' \end{pmatrix},$$

*where $\mathcal{Z}'$ is $\mathcal{N}(0, \Sigma^{(\mu)})$-distributed.*

(2) *If there exists $\gamma_1 \geq 0$ such that $\lim_{n\to\infty} a_n^{-1} c_n^6 = \gamma_1$, then*

$$\begin{pmatrix} \sqrt{a_n^{-1}c_n^2}(\theta_n - \theta) \\ \sqrt{\tilde{a}_n^{-1}}(\mu_n - \mu) \end{pmatrix} \xrightarrow{\mathcal{D}} \mathcal{N}\left( \begin{pmatrix} \sqrt{\gamma_1}\Delta^{(\theta)} \\ 0 \end{pmatrix}, \begin{pmatrix} \Sigma^{(\theta)} & 0 \\ 0 & \Sigma^{(\mu)} \end{pmatrix} \right).$$

*Comments on Theorem* 2. Set $(a_n) \equiv (a_0 n^{-1})$, $a_0$ satisfying (9), $(c_n) \equiv (c_0 n^{-1/6})$, $c_0 > 0$, and $(\tilde{a}_n) \equiv (\tilde{a}_0 n^{-1})$, $\tilde{a}_0 > 1/2$. Part (2) of Theorem 2 ensures that

$$\begin{pmatrix} n^{1/3}(\theta_n - \theta) \\ \sqrt{n}(\mu_n - \mu) \end{pmatrix} \xrightarrow{\mathcal{D}} \mathcal{N}\left( \begin{pmatrix} c_0^2\Delta^{(\theta)} \\ 0 \end{pmatrix}, \begin{pmatrix} a_0 c_0^{-2}\Sigma^{(\theta)} & 0 \\ 0 & \tilde{a}_0\Sigma^{(\mu)} \end{pmatrix} \right).$$

For this choice, $\theta_n$ converges with its optimal rate $n^{1/3}$, and $\mu_n$ converges with the parametric rate $\sqrt{n}$. Moreover, let us note that the asymptotic variance $\tilde{a}_0\Sigma^{(\mu)} = \tilde{a}_0^2[2\tilde{a}_0 - 1]^{-1}\sigma^2\delta^{-1}$ reaches its minimum $\sigma^2/\delta$ for $\tilde{a}_0 = 1$; the algorithm (4) is thus asymptotically efficient when $(\tilde{a}_n) \equiv (n^{-1})$.

To state the joint asymptotic behavior of $\overline{\theta}_n$ and $\mu_n$ defined in (5) and (7), we need to introduce the notation

$$R^{(\theta)} = \frac{1}{6}\left\{ \frac{\partial^3 f}{\partial x_i^3}(\theta) \right\}_{1 \leq i \leq d}, \tag{14}$$

as well as the following additional assumption.

(A5) (i) $\lim_{n\to\infty} \dfrac{na_n}{\log(\sum_{k=1}^n a_k)} = \infty$,

   (ii) $\lim_{n\to\infty} \dfrac{\sum_{k=1}^n a_k \log(\sum_{j=1}^k a_j)}{\sqrt{\sum_{k=1}^n c_k^2}} = 0$,

   (iii) $\lim_{n\to\infty} na_n^2 c_n^{-6} = \infty$.

THEOREM 3.  *Let $(\mu_n)$ be defined by* (7), *and assume that* (A1)–(A5) *hold with $(\tilde{a}_n) \equiv (n^{-1})$.*

(1) *If $\lim_{n\to\infty} nc_n^6 = \infty$, then*

$$\begin{pmatrix} c_n^{-2}(\overline{\theta}_n - \theta) \\ \sqrt{n}(\mu_n - \mu) \end{pmatrix} \xrightarrow{\mathcal{D}} \begin{pmatrix} -\left(\dfrac{1-2\tau}{1-4\tau}\right)[D^2 f(\theta)]^{-1} R^{(\theta)} \\ \mathcal{Z}' \end{pmatrix},$$

*where $\mathcal{Z}'$ is $\mathcal{N}(0, \sigma^2/\delta)$-distributed.*

(2) *If $\lim_{n\to\infty} nc_n^6 = 0$, then*

$$\begin{pmatrix} \sqrt{nc_n^2}(\overline{\theta}_n - \theta) \\ \sqrt{n}(\mu_n - \mu) \end{pmatrix} \xrightarrow{\mathcal{D}} \mathcal{N}\left( 0, \begin{pmatrix} \dfrac{(1-2\tau)\sigma^2}{2}[D^2 f(\theta)]^{-2} & 0 \\ 0 & \sigma^2/\delta \end{pmatrix} \right).$$

(3) *If there exists $\gamma_1 > 0$ such that $\lim_{n \to \infty} n c_n^6 = \gamma_1$, then*

$$\begin{pmatrix} \sqrt{n c_n^2}(\overline{\theta}_n - \theta) \\ \sqrt{n}(\mu_n - \mu) \end{pmatrix}$$

$$\xrightarrow{\mathcal{D}} \mathcal{N}\left( \begin{pmatrix} -2\gamma_1^{1/3}[D^2 f(\theta)]^{-1} R^{(\theta)} \\ 0 \end{pmatrix}, \begin{pmatrix} \dfrac{\gamma_1^{1/3}\sigma^2}{3}[D^2 f(\theta)]^{-2} & 0 \\ 0 & \sigma^2/\delta \end{pmatrix} \right).$$

Part (3) of Theorem 3 corresponds to the case where both $\overline{\theta}_n$ and $\mu_n$ are asymptotically efficient: they converge with their respective optimal rates $n^{1/3}$ and $n^{1/2}$, and their asymptotic covariance matrix is optimal (see, e.g., [8] for the optimality of the asymptotic covariance matrix of $\overline{\theta}_n$). To obtain the result of the third part of Theorem 3, one must choose $(c_n) \equiv (c_0 n^{-1/6})$, $c_0 > 0$, whereas different choices of the step size $(a_n)$ are possible. For instance, one may choose:

- either $(a_n) \equiv (a_0 n^{-\alpha})$, $a_0 > 0$, $\alpha \in \,]5/6, 1[$,
- or $(a_n) \equiv (a_0 n^{-1}[\log n]^\alpha)$, $a_0 > 0$, $\alpha > 0$,
- or $(a_n) \equiv (a_0 n^{-1}[\log \log n]^\alpha)$, $a_0 > 0$, $\alpha > 1$.

To conclude this section, let us mention that, in the case no additional observations are made to approximate $\mu$, averaging the algorithm (1) reduces the optimal convergence rate of the sequence $(\mu_n)$ then defined by (2). As a matter of fact, to average $\theta_n$, the step size $(a_n)$ in (1) must be chosen such that

$$\lim_{n \to \infty} n a_n / \log_2 n = \infty \tag{15}$$

[see assumption (A5)]. If the step size $(\tilde{a}_n)$ in (2) is set equal to $(n^{-1})$, then the combination of (A4) and (15) induces the condition $\lim_{n \to \infty} \tilde{a}_n^{-1} c_n^4 = \infty$, so that, in view of Theorem 1, $c_n^{-2}(\mu_n - \mu)$ converges to a degenerate distribution. Moreover, in this case the convergence rate $(c_n^{-2})$ is necessarily less than $\sqrt{n/(\log_2 n)^2}$. On the other hand, it is possible to choose $(\tilde{a}_n)$ such that $\tilde{a}_n^{-1/2}(\mu_n - \mu)$ converges to a Gaussian distribution. But, in this case also, because of the combination of (A4) and (15), the convergence rate $(\tilde{a}_n^{-1/2})$ is necessarily less than $\sqrt{n/(\log_2 n)^2}$. So, if the only parameter of interest in the double algorithm (1) and (2) is $\mu$, it is preferable not to average $\theta_n$: choosing in (1) the step size $(a_n) \equiv (n^{-1} \log_p n)$ (with $p > 2$) introduced by Koval and Schwabe [16] allows one to get rid of the tedious condition (9) on $a_0$ and to obtain better convergence rates for $(\mu_n)$ than those which can be achieved by averaging $\theta_n$.

**3. Proofs.** Let us first state some elementary properties of the classes $\mathcal{GS}(\alpha)$ of sequences that will be used throughout the proofs.

- If $(u_n) \in \mathcal{GS}(\alpha)$ and $(v_n) \in \mathcal{GS}(\beta)$, then $(u_n v_n) \in \mathcal{GS}(\alpha + \beta)$.
- If $(u_n) \in \mathcal{GS}(\alpha)$, then for all $c \in \mathbb{R}$, $(u_n^c) \in \mathcal{GS}(c\alpha)$.
- If $(u_n) \in \mathcal{GS}(\alpha)$, then for all $\epsilon > 0$ and $n$ large enough, $n^{\alpha - \epsilon} \leq u_n \leq n^{\alpha + \epsilon}$.
- If $(u_n) \in \mathcal{GS}(\alpha)$ and $\sum u_n = \infty$, then $\lim_{n \to \infty} n u_n [\sum_{k=1}^{n} u_k]^{-1} = 1 + \alpha$.

Now, set

$$(16) \qquad R_{n+1}^{(\theta)} = \frac{1}{2c_n} \{ f(\theta_n + c_n e_i) - f(\theta_n - c_n e_i) \}_{1 \leq i \leq d} - \nabla f(\theta_n),$$

$$(17) \qquad R_{n+1}^{(\mu)} = \begin{cases} \dfrac{1}{\delta} \sum_{i \in \mathcal{S}} [f(\theta_n + c_n e_i) + f(\theta_n - c_n e_i)] - \mu, \\ \qquad\qquad\qquad\qquad \text{if } (\mu_n) \text{ is defined by (2)}, \\ f(\theta_n) - \mu, \quad \text{if } (\mu_n) \text{ is defined by (4)}. \end{cases}$$

$$(18) \qquad \epsilon_{n+1}^{(\theta)} = \tfrac{1}{2} \{ W_{n,i}^+ - W_{n,i}^- \}_{1 \leq i \leq d} \in \mathbb{R}^d$$

and

$$(19) \qquad \epsilon_{n+1}^{(\mu)} = \begin{cases} \dfrac{1}{\delta} \sum_{i \in \mathcal{S}} [W_{n,i}^+ + W_{n,i}^-], \qquad \text{if } (\mu_n) \text{ is defined by (2)}, \\ \dfrac{1}{\delta} \sum_{i=1}^{\delta} \mathcal{W}_{n,i}, \qquad\qquad \text{if } (\mu_n) \text{ is defined by (4)}. \end{cases}$$

The recursive equation (1) can then be rewritten as

$$(20) \qquad \theta_{n+1} = \theta_n + a_n [\nabla f(\theta_n) + R_{n+1}^{(\theta)}] + \frac{a_n}{c_n} \epsilon_{n+1}^{(\theta)},$$

and the algorithms (2) and (4) as

$$(21) \qquad \mu_{n+1} = \mu_n + \tilde{a}_n [(\mu - \mu_n) + R_{n+1}^{(\mu)}] + \tilde{a}_n \epsilon_{n+1}^{(\mu)}.$$

These equations (20) and (21) can be viewed as particular stochastic approximation algorithms used for the search of a zero of a given function [of the function $\nabla f$ for (20) and of the function $x \mapsto \mu - x$ for (21)]. In Section 3.1, we state some preliminary results on stochastic approximation algorithms used for the search of zeros of a function $h$ that will be applied several times in the sequel; the proof of these preliminary results can be found in the technical report arxiv:math.ST/0610487v1. In Section 3.2 we establish an upper bound on the almost sure convergence rate of $\theta_n$, which will in particular be used to prove the strong consistency of $\mu_n$. In Section 3.3 we first prove Proposition 1, and then give an upper bound on the almost sure convergence rate of $\mu_n$ defined either by (2) or by (4). Section 3.4 is devoted to the proof of Theorems 1 and 2, and Section 3.5 to the proof of Theorem 3.

3.1. *Some preliminary results on stochastic approximation algorithms.*
We consider the stochastic approximation algorithm

$$Z_{n+1} = Z_n + \gamma_n[h(Z_n) + r_{n+1}] + \sigma_n \epsilon_{n+1}, \tag{22}$$

where the random variables $Z_0$, $(r_n)_{n\geq 1}$ and $(\epsilon_n)_{n\geq 1}$ are defined on a probability space $(\Omega, \mathcal{A}, \mathbb{P})$ equipped with a filtration $\mathcal{F} = (\mathcal{F}_n)$, and the step sizes $(\gamma_n)$ and $(\sigma_n)$ are two positive and nonrandom sequences that go to zero.

Stochastic approximation algorithms [such as (22)] used for the search of zeros of a function $h \colon \mathbb{R}^d \to \mathbb{R}^d$ have been widely studied under various assumptions; see [9, 23, 25] and the references therein. The object of this section is not to give the most general existing result on (22), but only to precisely state the results we shall use in the sequel for the study of (20) and (21); in particular, the hypotheses below are not the most general ones, but are appropriate in our framework.

(H1) There exists $z^* \in \mathbb{R}^d$ such that $\lim_{n\to\infty} Z_n = z^*$ a.s.

(H2) $h$ is differentiable at $z^*$, its Jacobian matrix $H$ at $z^*$ is symmetric, negative definite with maximal eigenvalue $-L < 0$, and there exists a neighborhood of $z^*$ in which $h(z) = H(z - z^*) + O(\|z - z^*\|^2)$.

(H3) (i) $\mathbb{E}(\epsilon_{n+1}|\mathcal{F}_n) = 0$ and there exists $m > 2$ such that $\sup_{n\geq 0} \mathbb{E}(\|\epsilon_{n+1}\|^m|\mathcal{F}_n) < \infty$.

   (ii) There exists a nonrandom, positive definite matrix $\Gamma$ such that $\lim_{n\to\infty} \mathbb{E}(\epsilon_{n+1}\epsilon_{n+1}^T|\mathcal{F}_n) = \Gamma$ a.s.

(H4) $r_{n+1} = R_{n+1}^{(1)} + O(\|Z_n - z^*\|^2)$ a.s., and there exist $\rho \in \mathbb{R}^d$ and a nonrandom sequence $(u_n)$ such that:

   (i) $\lim_{n\to\infty} \sqrt{u_n} R_{n+1}^{(1)} = \rho$ a.s.

   (ii) There exists $u^* > 0$ such that $(u_n) \in \mathcal{GS}(u^*)$.

(H5) (i) There exist $\alpha \in ]\max\{1/2, 2/m\}, 1]$ and $\beta > \alpha/2$ such that $(\gamma_n) \in \mathcal{GS}(-\alpha)$ and $(\sigma_n) \in \mathcal{GS}(-\beta)$.

   (ii) $\lim_{n\to\infty} n\gamma_n \in ]\max\{\frac{2\beta-\alpha}{2L}, \frac{u^*}{2L}\}, \infty]$, where $L$ and $u^*$ are defined in (H2) and (H4)(ii), respectively.

The asymptotic behavior of the algorithm (22) is given by the behavior of the sequences $(L_n)$ and $(\Delta_n)$ defined by

$$L_{n+1} = e^{(\sum_{k=1}^n \gamma_k)H} \sum_{k=1}^n e^{-(\sum_{j=1}^k \gamma_j)H} \sigma_k \epsilon_{k+1},$$

$$\Delta_{n+1} = (Z_{n+1} - z^*) - L_{n+1}.$$

In order to prove Proposition 1 and Theorems 1 and 2, we shall apply several times the following two lemmas.

LEMMA 1 [A.s. upper bound of $(L_n)$].   *Under hypotheses* (H2), (H3) *and* (H5), *we have* $\|L_n\| = O(\sqrt{\gamma_n^{-1}\sigma_n^2 \log(\sum_{k=1}^n \gamma_k)})$ *a.s.*

LEMMA 2 [A.s. convergence rate of $(\Delta_n)$].   *Under hypotheses* (H1)–(H5), *we have* $\lim_{n\to\infty} \sqrt{u_n}\Delta_n = -[H + \frac{\tilde{\xi}}{2}I_d]^{-1}\rho$ *a.s.*

Let us mention that, in particular, the combination of Lemmas 1 and 2 gives straightforwardly the following upper bound of the a.s. convergence rate of $Z_n$ toward $z^*$:

$$(23) \qquad \|Z_n - z^*\| = O\left(\sqrt{\gamma_n^{-1}\sigma_n^2 \log\left(\sum_{k=1}^n \gamma_k\right)} + u_n^{-1/2}\right) \qquad \text{a.s.}$$

To end this section, we now state a result concerning the averaged stochastic approximation algorithm derived from (22); we set

$$\overline{Z}_n = \frac{1}{\sum_{k=1}^n \gamma_k^2 \sigma_k^{-2}} \sum_{k=1}^n \gamma_k^2 \sigma_k^{-2} Z_k$$

and assume the following additional condition holds:

(H6)    (i)  $\lim_{n\to\infty} \dfrac{n\gamma_n}{\log(\sum_{k=1}^n \gamma_k)} = \infty.$

    (ii)  $\lim_{n\to\infty} \dfrac{\sum_{k=1}^n \gamma_k \log(\sum_{j=1}^k \gamma_j)}{\sqrt{\sum_{k=1}^n \gamma_k^2 \sigma_k^{-2}}} = 0.$

   (iii)  The sequence $(u_n)$ defined in assumption (H4) satisfies

$$\lim_{n\to\infty} n u_n \sigma_n^2 = \infty, \qquad \lim_{n\to\infty} \frac{\sum_{k=1}^n \gamma_k^2 \sigma_k^{-2} u_k^{-1}}{\sqrt{\sum_{k=1}^n \gamma_k^2 \sigma_k^{-2}}} = 0.$$

The asymptotic behavior of $(\overline{Z}_n)$ is given by the behavior of the sequences $(\Lambda_n)$ and $(\Xi_n)$ defined by

$$\Lambda_{n+1} = -\frac{1}{\sum_{k=1}^n \gamma_k^2 \sigma_k^{-2}} H^{-1} \sum_{k=1}^n \gamma_k \sigma_k^{-1} \epsilon_{k+1},$$

$$\Xi_{n+1} = (\overline{Z}_n - z^*) - \Lambda_{n+1}.$$

In Section 3.5, we shall apply several times the following lemma, which gives the asymptotic almost sure behavior of $(\Xi_n)$.

LEMMA 3 [A.s. convergence rate of $(\Xi_n)$].   *Assume that* (H1)–(H6) *hold.*
(1) *If* $\lim_{n\to\infty}[n\gamma_n^2\sigma_n^{-2}]^{-1/2}[\sum_{k=1}^n \gamma_k^2 \sigma_k^{-2} u_k^{-1/2}] = 0$, *then*

$$\lim_{n\to\infty} \sqrt{n\gamma_n^2\sigma_n^{-2}}\,\Xi_n = 0 \qquad a.s.$$

(2) *If the sequence $([n\gamma_n^2\sigma_n^{-2}]^{1/2}[\sum_{k=1}^n \gamma_k^2\sigma_k^{-2}u_k^{-1/2}]^{-1})$ is bounded, then*

$$\lim_{n\to\infty} \frac{n\gamma_n^2\sigma_n^{-2}}{\sum_{k=1}^n \gamma_k^2\sigma_k^{-2}u_k^{-1/2}}\ \Xi_n = -(1-2\alpha+2\beta)H^{-1}\rho \qquad a.s.$$

3.2. *Upper bound of the a.s. convergence rate of $\theta_n$.* Set

$$s_n = \sum_{k=1}^n a_k, \tag{24}$$

$$G = D^2 f(\theta), \tag{25}$$

$$L_{n+1}^{(\theta)} = e^{s_n G}\sum_{k=1}^n \frac{a_k}{c_k}e^{-s_k G}\epsilon_{k+1}^{(\theta)}, \tag{26}$$

$$\Delta_{n+1}^{(\theta)} = (\theta_{n+1}-\theta) - L_{n+1}^{(\theta)}. \tag{27}$$

The application of Lemma 1 to the recursive equation (20) [with $h \equiv \nabla f$, $(\gamma_n) \equiv (a_n)$ and $(\sigma_n) \equiv (a_n c_n^{-1})$] gives straightforwardly the following lemma.

LEMMA 4 [A.s. upper bound of $(L_n^{(\theta)})$]. *Under assumptions* (A2)(ii), (A3) *and* (A4)(i)–(iii), *we have* $\|L_n^{(\theta)}\| = O(\sqrt{a_n c_n^{-2}\log s_n})$ *a.s.*

Now, let $R_{n+1,i}^{(\theta)}$ denote the $i$th coordinate of $R_{n+1}^{(\theta)}$ [defined in (16)]; we have

$$\begin{aligned}
R_{n+1,i}^{(\theta)} &= \frac{1}{2c_n}\{[f(\theta_n+c_ne_i)-f(\theta_n)]-[f(\theta_n-c_ne_i)-f(\theta_n)]\} - \frac{\partial f}{\partial x_i}(\theta_n)\\
&= \frac{1}{2c_n}\bigg\{\bigg[c_n\frac{\partial f}{\partial x_i}(\theta_n)+\frac{c_n^2}{2}\frac{\partial^2 f}{\partial x_i^2}(\theta_n)+\frac{c_n^3}{6}\frac{\partial^3 f}{\partial x_i^3}(\theta_n)+o(c_n^3)\bigg]\\
&\quad -\bigg[-c_n\frac{\partial f}{\partial x_i}(\theta_n)+\frac{c_n^2}{2}\frac{\partial^2 f}{\partial x_i^2}(\theta_n)-\frac{c_n^3}{6}\frac{\partial^3 f}{\partial x_i^3}(\theta_n)+o(c_n^3)\bigg]\bigg\} - \frac{\partial f}{\partial x_i}(\theta_n)\\
&= \frac{c_n^2}{6}\frac{\partial^3 f}{\partial x_i^3}(\theta_n)+o(c_n^2),
\end{aligned}$$

and thus, in view of assumptions (A1) and (A2)(i), $\lim_{n\to\infty} c_n^{-2}R_{n+1}^{(\theta)} = R^{(\theta)}$ a.s., where $R^{(\theta)}$ is defined in (14). The application of Lemma 2 [with $(\sqrt{u_n}) \equiv (c_n^{-2})$ and $\rho \equiv R^{(\theta)}$] then gives the following lemma.

LEMMA 5 [A.s. convergence rate of $(\Delta_n^{(\theta)})$]. *Under assumptions* (A1)–(A3) *and* (A4)(i)–(iii), *we have* $\lim_{n\to\infty} c_n^{-2}\Delta_n^{(\theta)} = \Delta^{(\theta)}$ *a.s. where* $\Delta^{(\theta)}$ *is defined in* (11).

Let us note that the combination of Lemmas 4 and 5 ensures that, under assumptions (A1)–(A3) and (A4)(i)–(iii),

$$(28) \qquad \|\theta_n - \theta\| = O(\sqrt{a_n c_n^{-2} \log s_n} + c_n^2) \qquad \text{a.s.}$$

3.3. *On the a.s. asymptotic behavior of $\mu_n$ defined by* (2) *or* (4). In the case $\mu_n$ is defined either by (2) or by (4), the a.s. convergence of $\mu_n$ (resp. the a.s. convergence rate of $\mu_n$) is obtained by applying the Robbins–Monro theorem (resp. Lemmas 1 and 2) to the recursive equation (21). Since the $R_{n+1}^{(\mu)}$ term in (21) depends on $\theta_n$ [see (17)], we first upper bound this perturbation term by using the results of the previous section. To this end, we first note that in the case $(\mu_n)$ is defined by (2), we have

$$R_{n+1}^{(\mu)} = \frac{1}{\delta} \sum_{i \in \mathcal{S}} \left\{ \left[ f(\theta_n) + c_n \frac{\partial f}{\partial x_i}(\theta_n) + \frac{c_n^2}{2} \frac{\partial^2 f}{\partial x_i^2}(\theta_n) + o(c_n^2) \right] \right.$$

$$\left. + \left[ f(\theta_n) - c_n \frac{\partial f}{\partial x_i}(\theta_n) + \frac{c_n^2}{2} \frac{\partial^2 f}{\partial x_i^2}(\theta_n) + o(c_n^2) \right] \right\} - \mu$$

$$(29) \qquad = \frac{c_n^2}{\delta} \sum_{i \in \mathcal{S}} \frac{\partial^2 f}{\partial x_i^2}(\theta_n) + o(c_n^2) + [f(\theta_n) - f(\theta)]$$

$$= \frac{c_n^2}{\delta} \sum_{i \in \mathcal{S}} \frac{\partial^2 f}{\partial x_i^2}(\theta_n) + o(c_n^2) + O(\|\theta_n - \theta\|^2)$$

$$= \frac{c_n^2}{\delta} \sum_{i \in \mathcal{S}} \frac{\partial^2 f}{\partial x_i^2}(\theta_n) + o(c_n^2) + O\left( \frac{a_n \log s_n}{c_n^2} \right) \qquad \text{a.s.}$$

[where the last equality follows from the application of (28)]; in the case $(\mu_n)$ is defined by (4), similar computations give

$$(30) \qquad R_{n+1}^{(\mu)} = O\left( \frac{a_n \log s_n}{c_n^2} + c_n^4 \right) \qquad \text{a.s.}$$

In view of assumption (A4)(v), we deduce that:

- If $\lim_{n \to \infty} \tilde{a}_n^{-1} b_n^4 = 0$, then

$$(31) \qquad \lim_{n \to \infty} \sqrt{\tilde{a}_n^{-1}} R_{n+1}^{(\mu)} = 0 \qquad \text{a.s.}$$

- If $\lim_{n \to \infty} \tilde{a}_n^{-1} b_n^4 \in ]0, \infty]$, then

$$(32) \qquad \lim_{n \to \infty} \frac{1}{b_n^2} R_{n+1}^{(\mu)} = \frac{1}{\delta} \sum_{i \in \mathcal{S}} \frac{\partial^2 f}{\partial x_i^2}(\theta) \qquad \text{a.s.}$$

We can now prove Proposition 1 and give an upper bound on the a.s. convergence rate of $\mu_n$.

### 3.3.1. *Proof of Proposition* 1.

- In the case $\lim_{n\to\infty} \tilde{a}_n^{-1} b_n^4 = 0$, we have, in view of (31), $\tilde{a}_n |R_{n+1}^{(\mu)}|^2 = O(\tilde{a}_n^2)$ a.s., and thus, in view of (A4)(iv), $\sum \tilde{a}_n |R_{n+1}^{(\mu)}|^2 < \infty$ a.s.

- In the case $\lim_{n\to\infty} \tilde{a}_n^{-1} b_n^4 \in ]0, \infty]$, we have, in view of (32), $\tilde{a}_n |R_{n+1}^{(\mu)}|^2 = O(\tilde{a}_n b_n^4)$ a.s., and thus, in view of (A4)(v), $\sum \tilde{a}_n |R_{n+1}^{(\mu)}|^2 < \infty$ a.s.

In both cases, the application of the Robbins–Monro theorem (see, e.g., [9], page 61) ensures that $\sum \tilde{a}_n (\mu_n - \mu)^2 < \infty$ a.s. Since $\sum \tilde{a}_n = \infty$ [see (A4)(vi)], it follows that $\lim_{n\to\infty} \mu_n = \mu$ a.s.

### 3.3.2. *Upper bound on the a.s. convergence rate of $\mu_n$ defined by* (2) *or* (4). Set

$$(33) \qquad \tilde{s}_n = \sum_{k=1}^{n} \tilde{a}_k,$$

$$(34) \qquad L_{n+1}^{(\mu)} = e^{-\tilde{s}_n} \sum_{k=1}^{n} e^{\tilde{s}_k} \tilde{a}_k \epsilon_{k+1}^{(\mu)},$$

$$(35) \qquad \Delta_{n+1}^{(\mu)} = (\mu_{n+1} - \mu) - L_{n+1}^{(\mu)}$$

[where $\varepsilon_n^{(\mu)}$ is defined in (19)]. The application of Lemma 1 to the recursive equation (21) [with $h : x \mapsto \mu - x$, $(\gamma_n) \equiv (\tilde{a}_n)$ and $(\sigma_n) \equiv (\tilde{a}_n)$] gives straightforwardly the following lemma.

LEMMA 6 [A.s. upper bound of $(L_n^{(\mu)})$]. *Under assumptions* (A3), (A4)(iv) *and* (A4)(vi), *we have* $|L_n^{(\mu)}| = O(\sqrt{\tilde{a}_n \log \tilde{s}_n})$ *a.s.*

Moreover:

- if $\lim_{n\to\infty} \tilde{a}_n^{-1} b_n^4 = 0$, then, in view of (31), the application of Lemma 2 [with $(\sqrt{u_n}) \equiv (\sqrt{\tilde{a}_n^{-1}})$ and $\rho \equiv 0$] gives the first part of Lemma 7 below;
- if $\lim_{n\to\infty} \tilde{a}_n^{-1} b_n^4 \in ]0, \infty]$, then, in view of (32), the application of Lemma 2 [with $(\sqrt{u_n}) \equiv (b_n^{-2})$ and $\rho \equiv \frac{1}{\delta} \sum_{i \in \mathcal{S}} \frac{\partial^2 f}{\partial x_i^2}(\theta)$] gives the second part of Lemma 7 below.

LEMMA 7 [A.s. convergence rate of $(\Delta_n^{(\mu)})$]. *Let* (A1)–(A4) *hold.*

(1) *If* $\lim_{n\to\infty} \tilde{a}_n^{-1} b_n^4 = 0$, *then* $\lim_{n\to\infty} \sqrt{\tilde{a}_n^{-1}} \Delta_n^{(\mu)} = 0$ *a.s.*

(2) *If* $\lim_{n\to\infty} \tilde{a}_n^{-1} b_n^4 \in ]0, \infty]$, *then* $\lim_{n\to\infty} \frac{1}{b_n^2} \Delta_n^{(\mu)} = \Delta^{(\mu)}$ *a.s., where* $\Delta^{(\mu)}$ *is defined in* (13).

Although only Lemmas 6 and 7 will be used in the sequel, we state here the following proposition, which is obtained as a straightforward combination of these two lemmas, and which is of independent interest.

PROPOSITION 2 [A.s. upper bound of $(\mu_n - \mu)$].   *Under* (A1)–(A4), *we have*:
   (1) *If* $\lim_{n \to \infty} \tilde{a}_n^{-1} b_n^4 = 0$, *then* $|\mu_n - \mu| = O(\sqrt{\tilde{a}_n \log \tilde{s}_n})$ *a.s.*
   (2) *If* $\lim_{n \to \infty} \tilde{a}_n^{-1} b_n^4 \in ]0, \infty]$, *then* $|\mu_n - \mu| = O(\sqrt{\tilde{a}_n \log \tilde{s}_n} + b_n^2)$ *a.s.*

3.4. *Proof of Theorems* 1 *and* 2.   In view of the definition of $L_n^{(\theta)}$, $\Delta_n^{(\theta)}$, $L_n^{(\mu)}$ and $\Delta_n^{(\mu)}$ [see (26), (27), (34) and (35), resp.], Theorems 1 and 2 are straightforward consequences of the combination of Lemmas 5 and 7 together with the following lemma.

LEMMA 8 [Weak convergence rate of $(L_n^{(\theta)}, L_n^{(\mu)})$].   *Under* (A2)–(A4),

$$\begin{pmatrix} \sqrt{a_n^{-1} c_n^2} L_n^{(\theta)} \\ \sqrt{\tilde{a}_n^{-1}} L_n^{(\mu)} \end{pmatrix} \xrightarrow{\mathcal{D}} \mathcal{N}\left(0, \begin{pmatrix} \Sigma^{(\theta)} & 0 \\ 0 & \Sigma^{(\mu)} \end{pmatrix}\right),$$

*where* $\Sigma^{(\theta)}$ *and* $\Sigma^{(\mu)}$ *are defined in* (10) *and* (12), *respectively*.

PROOF.   Set

$$M_j^{(n)} = \begin{pmatrix} \sqrt{a_n^{-1} c_n^2} e^{s_n G} & 0 \\ 0 & \sqrt{\tilde{a}_n^{-1}} e^{-\tilde{s}_n} \end{pmatrix} \sum_{k=1}^{j} \begin{pmatrix} e^{-s_k G} a_k c_k^{-1} \epsilon_k^{(\theta)} \\ e^{\tilde{s}_k} \tilde{a}_k \epsilon_k^{(\mu)} \end{pmatrix}.$$

For each $n$, $M^{(n)} = (M_j^{(n)})_{j \geq 1}$ is a martingale whose predictable quadratic variation satisfies

$$\langle M \rangle_n^{(n)} = \begin{pmatrix} A_{1,n} & A_{2,n} \\ A_{2,n}^T & A_{4,n} \end{pmatrix}$$

with

$$A_{1,n} = a_n^{-1} c_n^2 e^{s_n G} \left\{ \sum_{k=1}^{n} a_k^2 c_k^{-2} e^{-s_k G} \mathbb{E}[\epsilon_k^{(\theta)} [\epsilon_k^{(\theta)}]^T | \mathcal{G}_{k-1}] e^{-s_k G^T} \right\} e^{s_n G^T},$$

$$A_{2,n} = a_n^{-1/2} \tilde{a}_n^{-1/2} c_n e^{s_n G} e^{-\tilde{s}_n} \left\{ \sum_{k=1}^{n} a_k \tilde{a}_k c_k^{-1} e^{-s_k G} e^{\tilde{s}_k} \mathbb{E}[\epsilon_k^{(\theta)} \epsilon_k^{(\mu)} | \mathcal{G}_{k-1}] \right\},$$

$$A_{4,n} = \tilde{a}_n^{-1} e^{-2\tilde{s}_n} \left\{ \sum_{k=1}^{n} \tilde{a}_k^2 e^{2\tilde{s}_k} \mathbb{E}[[\epsilon_k^{(\mu)}]^2 | \mathcal{G}_{k-1}] \right\}.$$

Now, under assumption (A3), we have, in view of (18) and (19),

$$\mathbb{E}[\epsilon_k^{(\theta)}[\epsilon_k^{(\theta)}]^T|\mathcal{G}_{k-1}] = \frac{\sigma^2 I_d}{2}, \qquad \mathbb{E}[[\epsilon_k^{(\theta)}]^T\epsilon_k^{(\mu)}|\mathcal{G}_{k-1}] = 0 \qquad \text{a.s.,}$$

$$\mathbb{E}[[\epsilon_k^{(\mu)}]^2|\mathcal{G}_{k-1}] = \frac{\sigma^2}{\delta} \qquad \text{a.s.}$$

It follows that $A_{2,n} = 0$ and, by application of Lemma 4 in [24], $\lim_{n\to\infty} A_{1,n} = \Sigma^{(\theta)}$ and $\lim_{n\to\infty} A_{4,n} = \Sigma^{(\mu)}$. We thus obtain

$$\lim_{n\to\infty} \langle M \rangle_n^{(n)} = \begin{pmatrix} \Sigma^{(\theta)} & 0 \\ 0 & \Sigma^{(\mu)} \end{pmatrix} \qquad \text{a.s.}$$

Moreover, in view of assumption (A3), we have

$$\sum_{k=1}^n \mathbb{E}[\|M_k^{(n)} - M_{k-1}^{(n)}\|^m|\mathcal{G}_{k-1}]$$

$$= O\left(\sum_{k=1}^n (a_n^{-1}c_n^2)^{m/2}\|e^{(s_n-s_k)G}a_k c_k^{-1}\|^m + \sum_{k=1}^n \tilde{a}_n^{-m/2}e^{-m(\tilde{s}_n-\tilde{s}_k)}\tilde{a}_k^m\right) \qquad \text{a.s.}$$

$$= O(w_n^{(\theta)} + w_n^{(\mu)}) \qquad \text{a.s.}$$

with

$$w_n^{(\theta)} = (a_n^{-1}c_n^2)^{m/2}e^{-mL^{(\theta)}s_n}\sum_{k=1}^n a_k^m c_k^{-m}e^{mL^{(\theta)}s_k},$$

$$w_n^{(\mu)} = \tilde{a}_n^{-m/2}e^{-m\tilde{s}_n}\sum_{k=1}^n e^{m\tilde{s}_k}\tilde{a}_k^m.$$

Now, since $(a_n^{-1}c_n^2) \in \mathcal{GS}(\alpha - 2\tau)$, we note that

$$w_{n+1}^{(\theta)} = \left[\frac{a_{n+1}^{-1}c_{n+1}^2}{a_n^{-1}c_n^2}\right]^{m/2}e^{-mL^{(\theta)}a_{n+1}}w_n^{(\theta)} + a_{n+1}^{m/2}$$

$$= \left[1 + \frac{\alpha - 2\tau}{n+1} + o\left(\frac{1}{n+1}\right)\right]^{m/2}[1 - mL^{(\theta)}a_{n+1} + o(a_{n+1})]w_n^{(\theta)} + a_{n+1}^{m/2}$$

$$= [1 + \xi^{(\theta)}a_{n+1} + o(a_{n+1})]^{m/2}[1 - mL^{(\theta)}a_{n+1} + o(a_{n+1})]w_n^{(\theta)} + a_{n+1}^{m/2}$$

$$= \left[1 + \frac{\xi^{(\theta)}m}{2}a_{n+1} + o(a_{n+1})\right][1 - mL^{(\theta)}a_{n+1} + o(a_{n+1})]w_n^{(\theta)} + a_{n+1}^{m/2}$$

$$= \left[1 - m\left(L^{(\theta)} - \frac{\xi^{(\theta)}}{2}\right)a_{n+1} + o(a_{n+1})\right]w_n^{(\theta)} + a_{n+1}^{m/2}.$$

Set $A^{(\theta)} \in ]0, L^{(\theta)} - \xi^{(\theta)}/2[$; for $n$ large enough, we get

$$|w_{n+1}^{(\theta)}| \leq (1 - A^{(\theta)}a_{n+1})|w_n^{(\theta)}| + a_{n+1}^{m/2},$$

and the application of Lemma 4.I.1 in [9] ensures that $\lim_{n\to\infty} w_n^{(\theta)} = 0$. In the same way, since $(\tilde{a}_n^{-1}) \in \mathcal{GS}(-\tilde{\alpha})$, we have

$$
\begin{aligned}
w_{n+1}^{(\mu)} &= \left[\frac{\tilde{a}_{n+1}}{\tilde{a}_n}\right]^{-m/2} e^{-m\tilde{a}_{n+1}} w_n^{(\mu)} + \tilde{a}_{n+1}^{m/2} \\
&= \left[1 - \frac{\tilde{\alpha}}{n+1} + o\left(\frac{1}{n+1}\right)\right]^{-m/2} [1 - m\tilde{a}_{n+1} + o(\tilde{a}_{n+1})] w_n^{(\mu)} + \tilde{a}_{n+1}^{m/2} \\
&= [1 - \xi^{(\mu)}\tilde{a}_{n+1} + o(\tilde{a}_{n+1})]^{-m/2} [1 - m\tilde{a}_{n+1} + o(\tilde{a}_{n+1})] w_n^{(\mu)} + \tilde{a}_{n+1}^{m/2} \\
&= \left[1 + \frac{\xi^{(\mu)}m}{2}\tilde{a}_{n+1} + o(\tilde{a}_{n+1})\right] [1 - m\tilde{a}_{n+1} + o(\tilde{a}_{n+1})] w_n^{(\mu)} + \tilde{a}_{n+1}^{m/2} \\
&= \left[1 - m\left(1 - \frac{\xi^{(\mu)}}{2}\right) a_{n+1} + o(a_{n+1})\right] w_n^{(\mu)} + a_{n+1}^{m/2},
\end{aligned}
$$

from which we deduce that $\lim_{n\to\infty} w_n^{(\mu)} = 0$. It thus follows that

$$
\sum_{k=1}^n \mathbb{E}[\|M_k^{(n)} - M_{k-1}^{(n)}\|^m | \mathcal{G}_{k-1}] = o(1) \qquad \text{a.s.,}
$$

and the application of Lyapounov's theorem gives

$$
M_n^{(n)} = \begin{pmatrix} \sqrt{a_n^{-1}c_n^2}\, L_n^{(\theta)} \\ \sqrt{\tilde{a}_n^{-1}}\, L_n^{(\mu)} \end{pmatrix} \xrightarrow{\mathcal{D}} \mathcal{N}\left(0, \begin{pmatrix} \Sigma^{(\theta)} & 0 \\ 0 & \Sigma^{(\mu)} \end{pmatrix}\right),
$$

which concludes the proof of Lemma 8.  □

3.5. *Proof of Theorem* 3. Set

$$
\Lambda_{n+1}^{(\theta)} = \frac{-1}{\sum_{k=1}^n c_k^2} G^{-1} \sum_{k=1}^n c_k \epsilon_{k+1}^{(\theta)},
$$

$$
\Xi_{n+1}^{(\theta)} = (\overline{\theta}_n - \theta) - \Lambda_{n+1}^{(\theta)},
$$

$$
\overline{\epsilon}_{k+1}^{(\mu)} = \frac{1}{\delta} \sum_{i=1}^{\delta} \overline{\mathcal{W}}_{k,i},
$$

$$
\Lambda_{n+1}^{(\mu)} = \frac{1}{n} \sum_{k=1}^n \overline{\epsilon}_{k+1}^{(\mu)},
$$

$$
\Xi_{n+1}^{(\mu)} = (\mu_{n+1} - \mu) - \Lambda_{n+1}^{(\mu)},
$$

where $\epsilon_k^{(\theta)}$ and $G$ are defined in (18) and (25), respectively. Theorem 3 follows straightforwardly from the combination of the three following lemmas, which give the a.s. convergence rate of $(\Xi_n^{(\theta)})$, of $(\Xi_n^{(\mu)})$ and the weak convergence rate of $(\Lambda_n^{(\theta)}, \Lambda_n^{(\mu)})$, respectively.

LEMMA 9 [A.s. convergence rate of $(\Xi_n^{(\theta)})$].  *Let the assumptions of Theorem 3 hold, and recall that $R^{(\theta)}$ is defined in (14).*

(1) *If $\lim_{n\to\infty} nc_n^6 = \infty$, then $\lim_{n\to\infty} c_n^{-2}\, \Xi_n^{(\theta)} = -(\frac{1-2\tau}{1-4\tau})G^{-1}R^{(\theta)}$ a.s.*

(2) *If $\lim_{n\to\infty} nc_n^6 = 0$, then $\lim_{n\to\infty} \sqrt{nc_n^2}\, \Xi_n^{(\theta)} = 0$ a.s.*

(3) *If there exists $\gamma > 0$ such that $\lim_{n\to\infty} nc_n^6 = \gamma$, then $\lim_{n\to\infty} \sqrt{nc_n^2}\, \Xi_n^{(\theta)} = -2\gamma^{1/3}G^{-1}R^{(\theta)}$ a.s.*

LEMMA 10 [A.s. convergence rate of $(\Xi_n^{(\mu)})$].  *Under the assumptions of Theorem 3 we have $\lim_{n\to\infty} \sqrt{n}\Xi_{n+1}^{(\mu)} = 0$ a.s.*

LEMMA 11 [Weak convergence rate of $(\Lambda_n^{(\theta)}, \Lambda_n^{(\mu)})$].  *Under the assumptions of Theorem 3, we have*

$$\begin{pmatrix} \sqrt{nc_n^2}\Lambda_n^{(\theta)} \\ \sqrt{n}\Lambda_n^{(\mu)} \end{pmatrix} \xrightarrow{\mathcal{D}} \mathcal{N}\left(0, \begin{pmatrix} \dfrac{\sigma^2(1-2\tau)}{2}G^{-2} & 0 \\ 0 & \dfrac{\sigma^2}{\delta} \end{pmatrix}\right).$$

PROOF OF LEMMA 9.  Set $(\gamma_n) \equiv (a_n)$, $(\sigma_n) \equiv (a_n c_n^{-1})$, $(u_n) \equiv (c_n^{-4})$ and $\epsilon \in ]0, (1-2\tau)/2[$. Since

$$\frac{\sum_{k=1}^{n} \gamma_k^2 \sigma_k^{-2} u_k^{-1}}{\sqrt{n\gamma_n^2 \sigma_n^{-2}}} = \frac{\sum_{k=1}^{n} c_k^6}{\sqrt{nc_n^2}} = O\left(\frac{n^\epsilon + nc_n^6}{\sqrt{nc_n^2}}\right) = o(1),$$

we can apply Lemma 3 to the recursive equation (20). Assumption (A4)(v) implies $\lim_{n\to\infty} nc_n^4 = \infty$, and thus $\sum c_n^4 = \infty$. Since $(c_n^4) \in \mathcal{GS}(-4\tau)$, we have

$$(36) \qquad \lim_{n\to\infty} \frac{nc_n^4}{\sum_{k=1}^{n} c_k^4} = 1 - 4\tau.$$

Consider the case $\lim_{n\to\infty} nc_n^6 \in ]0, \infty]$. We then have $\tau \leq 1/6$ and it follows from (36) that

$$\frac{\sqrt{n\gamma_n^2 \sigma_n^{-2}}}{\sum_{k=1}^{n} \gamma_k^2 \sigma_k^{-2} u_k^{-1/2}} = \frac{\sqrt{nc_n^2}}{\sum_{k=1}^{n} c_k^4} = O\left(\frac{1}{\sqrt{nc_n^6}}\right) = O(1).$$

The application of the second part of Lemma 3 then ensures that

$$\lim_{n\to\infty} \frac{nc_n^2}{\sum_{k=1}^{n} c_k^4}\Xi_n^{(\theta)} = -(1-2\tau)G^{-1}R^{(\theta)} \qquad \text{a.s.,}$$

and, applying (36) again, we obtain

$$(37) \qquad \lim_{n\to\infty} c_n^{-2}\Xi_n^{(\theta)} = -\left(\frac{1-2\tau}{1-4\tau}\right)G^{-1}R^{(\theta)} \qquad \text{a.s.,}$$

which gives the first part of Lemma 9. Note that if $\lim_{n \to \infty} nc_n^6 \in ]0, \infty[$, then $\tau = 1/6$; the third part of Lemma 9 follows straightforwardly from (37).

Now, consider the case $\lim_{n \to \infty} nc_n^6 = 0$. Set $\epsilon \in ]0, (1 - 2\tau)/2[$; using the fact that $(c_n^4) \in \mathcal{GS}(-4\tau)$ with $\tau \leq 1/4$ and applying (36) in the case $\tau \neq 1/4$, we obtain

$$\frac{\sum_{k=1}^{n} \gamma_k^2 \sigma_k^{-2} u_k^{-1/2}}{\sqrt{n \gamma_n^2 \sigma_n^{-2}}} = \frac{\sum_{k=1}^{n} c_k^4}{\sqrt{nc_n^2}} = O\left(\frac{n^\epsilon + nc_n^4}{\sqrt{nc_n^2}}\right) = o(1).$$

The application of the first part of Lemma 3 then ensures that $\lim_{n \to \infty} \sqrt{nc_n^2} \, \Xi_n^{(\theta)} = 0$ a.s., which concludes the proof of Lemma 9. $\square$

PROOF OF LEMMA 10. We have

$$|\Xi_{n+1}^{(\mu)}| = \left| \frac{1}{n} \sum_{k=1}^{n} \left[ \frac{1}{\delta} \sum_{i=1}^{\delta} Z_i(\overline{\theta}_k) \right] - f(\theta) - \Lambda_{n+1}^{(\mu)} \right|$$

$$= \left| \frac{1}{n} \sum_{k=1}^{n} [f(\overline{\theta}_k) - f(\theta)] \right|$$

$$= O\left( \frac{1}{n} \sum_{k=1}^{n} \|\overline{\theta}_k - \theta\|^2 \right)$$

$$= O\left( \frac{1}{n} \sum_{k=1}^{n} [\|\Lambda_{k+1}^{(\theta)}\|^2 + \|\Xi_{k+1}^{(\theta)}\|^2] \right).$$

By applying for instance Corollary 6.4.25 of [10], we get

$$\left\| \sum_{k=1}^{n} c_k \epsilon_{k+1}^{(\theta)} \right\| = O\left( \sqrt{\sum_{k=1}^{n} c_k^2 \log \log \left( \sum_{k=1}^{n} c_k^2 \right)} \right) \qquad \text{a.s.,}$$

and thus

$$\|\Lambda_{n+1}^{(\theta)}\|^2 = O((nc_n^2)^{-1} \log \log n) \qquad \text{a.s.}$$

The application of Lemma 9 then ensures that

$$|\Xi_{n+1}^{(\mu)}| = O\left( \frac{1}{n} \sum_{k=1}^{n} [c_k^4 + (kc_k^2)^{-1} \log \log k] \right) \qquad \text{a.s.}$$

$$= O(c_n^4 + (nc_n^2)^{-1} \log \log n) \qquad \text{a.s.}$$

In view of (A4)(v) (with $b_n = 0$ and $\tilde{a}_n = n^{-1}$), Lemma 10 follows. $\square$

Proof of Lemma 11. Set

$$\mathcal{M}_j^{(n)} = \begin{pmatrix} -\left(\sum_{k=1}^n c_k^2\right)^{-1/2} G^{-1} & 0 \\ 0 & n^{-1/2} \end{pmatrix} \sum_{k=1}^j \begin{pmatrix} c_k \epsilon_k^{(\theta)} \\ \overline{\epsilon}_k^{(\mu)} \end{pmatrix}.$$

In view of (A3), for each $n$, $\mathcal{M}^{(n)} = (\mathcal{M}_j^{(n)})_{j \geq 1}$ is a martingale whose predictable quadratic variation satisfies

$$\langle \mathcal{M} \rangle_n^{(n)} = \begin{pmatrix} \dfrac{\sigma^2}{2} G^{-2} & 0 \\ 0 & \dfrac{\sigma^2}{\delta} \end{pmatrix} \qquad \text{a.s.}$$

and we have

$$\sum_{k=1}^n \mathbb{E}[\|\mathcal{M}_k^{(n)} - \mathcal{M}_{k-1}^{(n)}\|^m | \mathcal{G}_{k-1}]$$

$$= O\left(\left[\sum_{k=1}^n c_k^2\right]^{-m/2} \sum_{k=1}^n c_k^m + n^{1-m/2}\right) \qquad \text{a.s.}$$

$$= o(1) \qquad \text{a.s.}$$

The application of Lyapounov's theorem then ensures that

$$\mathcal{M}_n^{(n)} = \begin{pmatrix} \sqrt{\sum_{k=1}^n c_k^2} \Lambda_n^{(\theta)} \\ \sqrt{n}\, \Lambda_n^{(\mu)} \end{pmatrix} \xrightarrow{\mathcal{D}} \mathcal{N}\left(0, \begin{pmatrix} \dfrac{\sigma^2}{2} G^{-2} & 0 \\ 0 & \dfrac{\sigma^2}{\delta} \end{pmatrix}\right),$$

and Lemma 11 follows from the fact that, since $(c_n^2) \in \mathcal{GS}(-2\tau)$ with $\tau > 1/2$, we have $\lim_{n \to \infty} n c_n^2 [\sum_{k=1}^n c_k^2]^{-1} = 1 - 2\tau$. $\square$

**Acknowledgments.** The authors deeply thank an Associate Editor and a referee for helpful comments and advice that led to a substantial improvement of the original version of this paper.

DÉPARTEMENT DE MATHÉMATIQUES, BAT. FERMAT
UNIVERSITÉ DE VERSAILLES–SAINT-QUENTIN
45, AVENUE DES ETATS-UNIS
78035 VERSAILLES CEDEX
FRANCE
E-MAIL: mokkadem@math.uvsq.fr
           pelletier@math.uvsq.fr